\documentclass[11pt,a4paper]{article}

\usepackage[utf8]{inputenc}
\usepackage[T1]{fontenc}
\usepackage[margin=1in]{geometry}

\usepackage{amsmath,amssymb,amsthm}
\usepackage{graphicx}
\usepackage{epstopdf}% To incorporate .eps illustrations using PDFLaTeX, etc.
\usepackage[caption=false]{subfig}% Support for small, `sub' figures and tables
\usepackage{booktabs}% used by the original class
\usepackage{authblk}% author / affiliation blocks
\usepackage{qubo}% your own macros

\usepackage[numbers,sort&compress]{natbib}
\bibpunct[, ]{[}{]}{,}{n}{,}{,}
\renewcommand\bibfont{\fontsize{10}{12}\selectfont}

\usepackage{titlesec}% to format the appendix headings

\hypersetup{colorlinks=true,allcolors=blue}

\theoremstyle{plain}

\theoremstyle{definition}

\theoremstyle{remark}

\title{\texttt{QUBO.jl}: A Julia Ecosystem for Quadratic Unconstrained
Binary Optimization\thanks{The published version is available at
\href{https://doi.org/10.1080/10556788.2026.2702926}%
{https://doi.org/10.1080/10556788.2026.2702926}.}}

\author[1,2,6]{Pedro Maciel Xavier\thanks{Pedro Maciel Xavier was
affiliated with institutions 1, 2, and 6 during the development of this
work.}}
\author[2,3]{Pedro Ripper}
\author[2]{Tiago Andrade}
\author[2]{Joaquim Dias Garcia}
\author[1]{Nelson Maculan}
\author[4,5,6]{David E. Bernal Neira\thanks{Corresponding author:
\href{mailto:dbernaln@purdue.edu}{dbernaln@purdue.edu}.}}

\affil[1]{PESC -- COPPE -- Federal University of Rio de Janeiro,
Av.\ Horácio Macedo, 2030, Rio de Janeiro, RJ 21941-598, Brazil}
\affil[2]{PSR, Praia de Botafogo, 370, Rio de Janeiro, RJ 22250-040, Brazil}
\affil[3]{DEE -- Pontifical Catholic University of Rio de Janeiro,
Rua Marquês de São Vicente, 225, Rio de Janeiro, RJ 22451-900, Brazil}
\affil[4]{Research Institute of Advanced Computer Science, Universities Space
Research Association, 425 3rd Street SW, Suite 950, Washington, DC 20024, USA}
\affil[5]{Quantum Artificial Intelligence Laboratory, NASA Ames Research
Center, PO Box 1, Moffett Field, CA 94035, USA}
\affil[6]{Davidson School of Chemical Engineering, Purdue University,
480 Stadium Mall Drive, West Lafayette, IN 47907, USA}

\date{}

\begin{document}

\maketitle

\begin{abstract}
We present \texttt{QUBO.jl}, an end-to-end Julia package for working with QUBO (Quadratic Unconstrained Binary Optimization) instances.
This tool aims to convert a broad range of optimization problems in \texttt{JuMP}, Julia's mathematical programming package, for straightforward application in many physics and physics-inspired solution methods whose standard model form is equivalent to QUBO.
These methods include quantum annealing, quantum gate-circuit optimization algorithms (Quantum Optimization Alternating Ansatz, Variational Quantum Eigensolver), other hardware-accelerated platforms, such as Coherent Ising Machines and Simulated Bifurcation Machines, and more traditional methods such as simulated annealing.
In addition to working with reformulations, \texttt{QUBO.jl} allows its users to interface with the aforementioned hardware, sending QUBO models to these devices and retrieving results for subsequent analysis.
\texttt{QUBO.jl} was written as a JuMP / MathOptInterface (MOI) layer that automatically maps between the input and output frames, thus providing a smooth modeling experience.
\end{abstract}


\section{Introduction}\label{sec:intro}
Over the past two decades, the mathematical programming framework Quadratic Unconstrained Binary Optimization (QUBO)~\cite{Kochenberger2014} has gained popularity and relevance in many research fields such as finance~\cite{qc_finance, herman2023quantum}, chemistry~\cite{cao2019quantum}, engineering~\cite{bernal_perspectives_2022}, and others~\cite{Biamonte_2017,rieffel2019ansatze}.
This trend of increasing attention arose primarily from significant advances in computing paradigms whose standard task is to sample high-quality solutions to this kind of problem.
In addition to that, QUBO formulations are known to be well-suited when representing many non-convex global optimization problems, especially those of combinatorial and discrete nature~\cite{qc_finance,rieffel2019ansatze,bernal_perspectives_2022, VERMA2022100594}.

In this regard, quantum computers stand out as one of the leading platforms for running QUBO formulations. 
Although not yet presenting real-world application results, quantum computing (QC) has been gaining a lot of traction, backed by investments from the public and private sectors interested in the expected benefits of using quantum computers~\cite{quantum_landscape}.

% From this point of view, the current state of QC, known as the Noisy Intermediate Quantum~(NISQ) Era~\cite{Preskill_2018} is characterized by big tech companies, startups, and national institutes that improve and develop quantum systems, which are currently susceptible to noisy operations and decoherence~\cite{quantum_landscape}.

However, in a few years, QC is expected to show some practical results in areas such as optimization~\cite{Preskill_2018}. 
Moreover, exciting proofs of concept have emerged, highlighting the potential of quantum technologies to address problems of this form~\cite{tasseff2022emerging}.
This advent will increase the demand for a trained workforce to program these quantum computers, due to the fact that Operations Research specialists usually lack the required background. 
As a consequence, the ability to adapt to this new scenario will be key to thriving with this new framework of computation, requiring companies and research institutes to spend time and resources staying on track with technology.

For this undertaking, one challenge is to reformulate optimization models into a format amenable to Quantum Computing.
Furthermore, Operations Research specialists must be capable of operating on these new architectures, where QUBO arises as their standard optimization format, although presenting different interfaces.
In addition, it will be crucial to transform the raw results of different hardware into helpful information to analyze.

% To address these issues, we have developed the \href{https://github.com/JuliaQUBO/QUBO.jl}{\texttt{QUBO.jl}} package~\cite{qubojl:2023}, a comprehensive tool for working with QUBO problems. 
% This software works as a wrapper for three other packages. The first, \href{https://github.com/JuliaQUBO/ToQUBO.jl}{\texttt{ToQUBO}}~\cite{toqubo:2022} handles the reformulation of general optimization problems in the QUBO form. 
% The second \href{https://github.com/JuliaQUBO/QUBODrivers.jl}{\texttt{QUBODrivers.jl}}~\cite{qubodrivers:2023} is responsible for hardware interfacing, sending QUBOs to quantum devices and analyzing generated solution pools.
% Finally, the third one, \href{https://github.com/JuliaQUBO/QUBOTools.jl}{\texttt{QUBOTools.jl}}~\cite{qubotools:2023} was envisioned to handle QUBO file conversions from different architectures and provide some utilities for the other two packages, such as graphing recipes.

To address these issues, we have developed the \href{https://github.com/JuliaQUBO/QUBO.jl}{\texttt{QUBO.jl}} package~\cite{qubojl:2023}, a comprehensive tool for working with QUBO problems. 
This software works as a wrapper for three other packages, \href{https://github.com/JuliaQUBO/ToQUBO.jl}{\texttt{ToQUBO}}~\cite{toqubo:2022}, {\texttt{QUBODrivers.jl}}~\cite{qubodrivers:2023} and {\texttt{QUBOTools.jl}}~\cite{qubotools:2023}, each performing a specific task.

Having all three packages combined into \texttt{QUBO.jl}, the user will be able to evaluate the potential of quantum and hardware-accelerated computers when addressing their problems without prior knowledge of these systems, providing a simple introduction to optimization via hardware accelerators, such as quantum computers.

In this work, we will present \texttt{QUBO.jl}'s packages separately. 
This paper is organized as follows.
Section \ref{sec:compilation} briefly discusses the adopted reformulation methods for generating QUBO problems from a more general class, such as mixed-integer nonlinear problems~(MINLP).
Section \ref{sec:virtual} follows with a general explanation of our reformulation layer, which is responsible for providing a smooth modeling experience to the user who is familiar with the \texttt{Julia} mathematical programming modeling package, \href{https://github.com/jump-dev/JuMP.jl}{\texttt{JuMP}}~\cite{Lubin2023}.
The proposed QUBO sampler interface for \texttt{JuMP}, which is implemented in the \texttt{QUBODrivers.jl} module, is presented next in Section \ref{sec:solver}.
Next, Section \ref{sec:format-tools} briefly introduces our package for handling different QUBO file formats, \texttt{QUBOTools.jl}.
The current state of quantum software packages is analyzed in Section \ref{sec:julia_env}, while Section \ref{sec:comparison} sets a benchmark with other existing projects, and Section \ref{sec:conclusions} points out conclusions and next steps in the development of this project.
% We also provide in Appendix \ref{sec:appendix:coefficients} a discussion on the comparison of integer to binary encoding, while our new encoding method is presented in Appendix \ref{sec:appendix:ap}.
Finally, we provide details on our benchmarking environment in Appendix \ref{sec:appendix:benchmark}.

\subsection{Quadratic Unconstrained Binary Optimization (QUBO)}\label{sec:qubo}

An optimization problem is in the QUBO form if it can be written as:
\begin{equation}
\begin{array}{>{\displaystyle}r>{\displaystyle}l}
(\textrm{QUBO}): \min_{\vec{x} \in \mathbb{B}^{n}} & \vec{x}' \vec{Q}\, \vec{x} = \min_{[x_1, \cdots, x_n] \in \mathbb{B}^n} 2 \sum_{i = 1}^{n} \sum_{j = i + 1}^{n} q_{i, j}\, x_i\, x_j + \sum_{i = 1}^{n} q_{i, i}\, x_i  
    \label{eq:qubo}
\end{array}
\end{equation}
where $\vec{Q} \in \mathbb{R}^{n \times n}$ is a symmetric matrix with elements $q_{i,j}$.
It is worth noting that $\vec{Q}$ is not required to be symmetric for a general QUBO representation. 
The main requirement is for the objective function to be at most quadratic, comprised of only binary variables, and be an unconstrained problem.
% When $\vec{Q}$ is not diagonal, the solution of this problem is known to be, in general, NP-Hard~\cite{pardalos1992complexity}.

\subsection{Ising Model}\label{sec:ising}

An optimization problem is in the Ising form if it is written as:
\begin{equation}
\begin{array}{>{\displaystyle}r>{\displaystyle}l}
(\textrm{Ising}): \min_{\vec{s} \in \{\pm 1\}^{n}} & \vec{s}' \vec{J}\, \vec{s} + \vec{h}' \vec{s} = \min_{[s_1, \cdots, s_n] \in  \{\pm 1\}^{n}} \sum_{i = 1}^{n} \sum_{j = i + 1}^{n} J_{i, j}\, s_i\, s_j + \sum_{i = 1}^{n} h_{i}\, s_i,
    \label{eq:ising}
\end{array}
\end{equation}
where $\vec{J} \in \mathbb{R}^{n \times n}$ is strictly upper triangular and $\vec{h} \in \mathbb{R}^{n}$.

There is a bijective relation between Ising and QUBO models, where one could easily encode variables from one model to the other as:
\begin{align*}
    \vec{s} &= 2 \vec{x} - \vec{1} \\
    \vec{x} &= \frac{1}{2} (\vec{s} + \vec{1})
\end{align*}

\subsection{Mixed-Integer Nonlinear Programming}\label{sec:miqp}

An optimization problem is in the Mixed-Integer Nonlinear Programming (MINLP) form if it can be written as:
\begin{equation}
\begin{array}{>{\displaystyle}r>{\displaystyle}l>{\displaystyle}l}
(\textrm{MINLP}): \min_{\vec{y}, \vec{z}} & f(\vec{y}, \vec{z}) \\
    \textrm{s.t.}                    & \vec{g}(\vec{y}, \vec{z}) \le \vec{0} \\
    ~                                & \vec{h}(\vec{y}, \vec{z}) = \vec{0}   \\
    ~                                &\vec{y} \in \mathbb{R}^{n}; \vec{z} \in \mathbb{Z}^{m},    \\
\end{array}
\end{equation}
where \(f: \mathbb{R}^{n} \times \mathbb{Z}^{m} \to \mathbb{R}\) denotes the objective function, which by convention is defined to be minimized although \( \displaystyle \min_{\vec{y}, \vec{z}} ~ f(\vec{y}, \vec{z}) = -\max_{\vec{y}, \vec{z}} ~ -f(\vec{y}, \vec{z}) \), \(\vec{y}\) and \(\vec{z}\) denote continuous and integer variables, respectively, and \(\vec{g}\) and \(\vec{h}\) denote inequality and equality vector constraints, respectively.
%Inequality constraints can always be posed as equalities by adding a non-negative variable to the left-hand side corresponding to the slackness of the constraint.

Notice that we do not make any assumption about the convexity or linearity of the functions \(f, \vec{g}, \vec{h}\), making MINLP a very flexible modeling paradigm.
The possibility of modeling linear and nonlinear constraints together with discrete and continuous variables makes it able to represent Universal Turing Machines~\cite{liberti2014mathematical}, and has many applications in science and engineering~\cite{belotti2013mixed,kronqvist2019review}.

{Considering the types of constraints that \texttt{QUBO.jl} provides an automatic reformulation for, which will be discussed in Section \ref{sec:comparison}, we are particularly interested in cases where the nonlinear objective and constraints can be written as quadratics, i.e., \(f(\vec{y},\vec{z}) = [\vec{y};\vec{z}]' \vec{Q} [\vec{y};\vec{z}] + \vec{\ell}'[\vec{y};\vec{z}] + c \).}
A problem with both quadratic constraints and quadratic objective is known as a Mixed-Integer Quadratically Constrained Quadratic Program (MIQCQP), while the case with linear constraints, i.e., \(\vec{h}(\vec{y},\vec{z}) = A[\vec{y};\vec{z}] - \vec{b}\), and with only a quadratic objective and linear constraints is known as a Mixed-Integer Quadratic Program (MIQP).
Thus, \texttt{QUBO.jl} currently provides automatic reformulation for MINLP models that fall within the subclass of MIQCQP and MIQP.
\begin{equation}
\begin{array}{>{\displaystyle}r>{\displaystyle}l>{\displaystyle}l}
(\textrm{MIQP}): \min_{\vec{y}, \vec{z}} & [\vec{y};\vec{z}]' \vec{Q}\, [\vec{y};\vec{z}] + \vec{\ell}'[\vec{y};\vec{z}] + c \\
    \textrm{s.t.}                    & A[\vec{y};\vec{z}]  = \vec{b} \\
    ~                                & \vec{y} \in \mathbb{R}^{n}; \vec{z} \in \mathbb{Z}^{m}   \\
\end{array}
\end{equation}
where \( \vec{Q} \in \mathbb{R}^{(n+m) \times (n+m)} \) is symmetric, \( \vec{\ell} \in \mathbb{R}^{(n+m)} \), \( c \in \mathbb{R} \), \( A \in \mathbb{R}^{(n+m) \times p} \), \( \vec{b} \in \mathbb{R}^{p} \).

MIQP models are notorious for encompassing both discrete and simple nonlinear behavior.
They are widely used in areas such as Finance, Economics, Computer Science, Operations Research, and Engineering~\cite{misener2012global}.
Every QUBO can also be seen as an MIQP since \( \mathbb{B} \subset \mathbb{Z} \).
This aspect allows one to use commercial MIQP solvers, such as Xpress~\cite{xpress} and CPLEX~\cite{cplex2009v12}, to work with QUBO instances.
It is worth mentioning that some optimization software, e.g., Gurobi~\cite{gurobi} and SCIP~\cite{rehfeldt2023faster}, have released specific routines to handle QUBO models.
On the other hand, it is possible to approximate MIQP problems~\cite{liberti2009reformulations} through QUBO modeling by applying specific variable encoding and constraint mapping methods, as discussed in Section \ref{sec:compilation}.

\section{\texttt{QUBO.jl} ecosystem}

\texttt{QUBO.jl} was developed to become a bridge between Operations Research specialists and quantum or quantum-inspired optimization platforms.
Wrapping the previously mentioned packages, \texttt{ToQUBO.jl}, \texttt{QUBODrivers.jl} and \texttt{QUBOTools.jl}, it turns into a self-contained ecosystem that can cover most tasks related to QUBO and Ising optimization.
From the diagram depicted in Figure \ref{fig:virtual-pipeline-simple}, one can better understand how \texttt{QUBO.jl} works and how its packages communicate.

First, \texttt{ToQUBO.jl} acts as an interface for \texttt{JuMP} users to generate QUBO instances.
As they send their \texttt{JuMP} model to \texttt{ToQUBO.jl}, after reformulating, it caches the QUBO problem as a \texttt{MOI} model~(a lower-level representation of an Optimization problem, based on \texttt{MOI}'s data structures).
Then, this \texttt{MOI} model can either be sent to a general classical MIQP solver or, as we have labeled, a QUBO solver that can be a classical, quantum, or quantum-inspired hardware.

{  In order for a QUBO model, generated by \texttt{ToQUBO.jl}, to be sent to a QUBO solver, we have developed \texttt{QUBODrivers.jl}, where users can define an interface to access different QUBO-amenable solvers, such as gate-based algorithms~\cite{cerezo2021variational}, Quantum Annealing~\cite{McGeoch2014}, and other quantum-inspired methods~\cite{mohseni2022ising}.}
After sending a QUBO to one of these solvers, \texttt{QUBODrivers.jl} retrieves the results for the user, which can be later evaluated by some of its analysis tools or sent to \texttt{QUBOTools.jl}, which provides plotting recipes for a visual representation of the results.

As some emerging architectures require QUBOs to be sent in files, another feature from \texttt{QUBOTools.jl} is file-conversion handling.
This package allows format conversion between different file formats, considering that each hardware demands a specific type.
\texttt{QUBOTools.jl} also provides translating \texttt{MOI} QUBO models into a file and vice-versa, presenting a bidirectional conversion between any of the envisioned QUBO model representations.

In summary, \texttt{QUBO.jl} users can harness quantum and quantum-inspired optimization methods to sample solutions for any general optimization problem. 
Therefore, we aim to provide a frictionless experience in integrating our package into the work of Operations Research specialists, requiring minimal previous knowledge of quantum technologies.

\begin{figure}[hb]
    \centering
    \begin{tikzpicture}[%
        model/.style = {%
                draw,
                thick,
                align=center,
                color=black,
                minimum width=15mm,
                minimum height=10mm,
                shape=rectangle,
                inner sep=2.5ex,
            },%
        solver/.style = {%
                fill           = white,%
                draw           = black,%
                very thick,%
                text centered,%
                rounded corners,
                % shape = diamond,%
                align          = center,%
            },
        arrow/.style = {%
            very thick,%
            % shorten <=0.025cm,%%
            % shorten >=0.025cm,%%
            preaction = {%
                draw, -, white,%%
                line width=0.25cm,%%
                shorten <=0.1cm,%%
                shorten >=0.1cm,%%
            }
        },%
        component/.style = {%
                very thick,%
                text centered,%
                fill           = white,%
                draw           = black,%
                minimum height = 2cm,%
                minimum width  = 2.5cm,%
                align          = center,%
            },
        legend/.style = {%
                thick,%
                fill           = white,%
                draw           = gray,%
                minimum height = 1cm,%
                minimum width  = 3.5cm,%
                align          = left,%
        },
        data/.style = {%
                very thick,%
                cylinder,%
                text centered,%
                fill                = white,%
                draw                = black,%
                shape border rotate = 90,%
                minimum height      = 1cm,%
                minimum width       = 2cm,%
                align               = center,%
            },
        misc/.style = {%
                very thick,%
                shape=ellipse,%
                text centered,%
                fill                = white,%
                draw                = black,%
                shape border rotate = 90,%
                minimum height      = 1cm,%
                minimum width       = 1cm,%
                align               = center,%
            }
    ]
    \footnotesize%
    \node[model] (JuMP) at (0, -3) {%
        \textbf{\texttt{JuMP} MINLP Model} \\[2ex]%
        \footnotesize%
        $\begin{array}{>{\displaystyle}r>{\displaystyle}l}
                \textrm{opt}  & f(\vec{y}, \vec{z})        \\
                \textrm{s.t.} & g(\vec{y}, \vec{z}) \le 0  \\
                ~             & h(\vec{y}, \vec{z}) = 0    \\
                ~             & \vec{y} \in \mathbb{R}^{m}; \vec{z} \in \mathbb{Z}^{n}
            \end{array}$%
    };
    \node[model] (MOI) at (0, -5.5) {%
        \textbf{\texttt{MOI} QUBO Model} \\[2ex]%
        \footnotesize%
        $\begin{array}{>{\displaystyle}r>{\displaystyle}l}
                \textrm{opt}  & \vec{x}' \vec{Q} \vec{x}         \\
                \textrm{s.t.} & \vec{x} \in \mathbb{B}^{n}
            \end{array}$%
    };

    \node[
        component,%
        ultra thick,%
        minimum height=7.5cm,%
        minimum width=7.5cm,%
    ] (QUBO) at (6, -5) {\normalsize\bfseries\texttt{QUBO.jl}};
    
    \node[component] (ToQUBO) at (4, -3) {%
        \texttt{ToQUBO.jl}%
    };

    \node[component] (QUBODrivers) at (4, -7) {%
        \texttt{QUBODrivers.jl}%
    };

    \node[component] (QUBOTools) at (8, -7) {%
        \texttt{QUBOTools.jl}%
    };

    \node[solver] (MIQP) at (0, -7.5) {%
        \textbf{MIQP Solver}%
    };

    \node[solver] (QUBOSolver) at (0, -8.5) {%
        \textbf{QUBO Solver}%
    };

    \node[data] (File) at (11.5, -5) {%
        \textbf{File}%
    };

    \node[misc] (Analysis) at (11.5, -7) {%
        \textbf{Analysis}%
    };

    % \draw[arrow, {Kite[open]}-{Kite}] (JuMP) -- (MOI);
    \draw[arrow, {Kite[open]}-{Kite}] (JuMP) -- (ToQUBO);
    \draw[arrow, {Kite[open]}-{Kite}] (ToQUBO.-90) |- (MOI.20);
    \draw[arrow, {Kite[open]}-{Kite}] (MOI) -- (MIQP);
    \draw[arrow, {Kite[open]}-{Kite}] (MOI) -| (QUBODrivers.90);
    \draw[arrow, {Kite[open]}-{Kite}] (QUBODrivers) |- (QUBOSolver);
    \draw[arrow, {Kite}-{Kite}] (File) -| (QUBOTools.60);
    \draw[arrow,  -{Kite}] (ToQUBO.20) -| (QUBOTools);
    \draw[arrow,  -{Kite[open]}] (ToQUBO.-20) -| (QUBOTools.120);
    \draw[arrow, {Kite}-{Kite[open]}] (QUBODrivers) -- (QUBOTools);
    \draw[arrow,  -{Kite[open]}] (QUBOTools) -- (Analysis);

    \node[legend] (Legend) at (12, -2.25){};

    \node[minimum width=2cm] (ModelArrow) at (12.5, -2.0) {Model Data};
    \draw[arrow, -{Kite}] (10.5, -2.0) -- (ModelArrow);
    
    \node[minimum width=2cm] (SolutionArrow) at (12.5, -2.5) {Solution Data};
    \draw[arrow, -{Kite[open]}] (10.5, -2.5) -- (SolutionArrow);

\end{tikzpicture}
    \caption{Diagram of the \texttt{QUBO.jl} ecosystem. \texttt{ToQUBO.jl} acts as an interface for QUBO reformulation within the \texttt{JuMP} framework. The QUBO model generated by \texttt{ToQUBO.jl}, or explicitly defined in a \texttt{JuMP} problem or a file, can later be sent to QUBO solvers with \texttt{QUBODrivers.jl}. Finally, \texttt{QUBOTools.jl} provides tools for result and model analysis. \\}
    \label{fig:virtual-pipeline-simple}
\end{figure}

\section{QUBO Compilation}\label{sec:compilation}
\texttt{ToQUBO.jl}'s main goal is to translate \texttt{JuMP} models into the QUBO form, aiming at their submission to solution sampling architectures, such as quantum annealers.
To illustrate the proposed pipeline, it might be worth building an analogy by pretending that a QUBO instance is equivalent to “assembly code” when dealing with some of the optimization machines presented above.
Within this picture, \texttt{ToQUBO.jl} acts as a higher-level language compiler, allowing more general \texttt{JuMP} models to be used as input for solvers that consider the QUBO formalism.
Almost as \texttt{gcc}~\cite{stallman1987c} generates Assembly code from programs written in \texttt{C}, \texttt{ToQUBO.jl} is capable of producing optimization programs as a form of lower code for specialized hardware.

Conceptually speaking, each characteristic of the target model introduces one or more specific steps to the reformulation process.
The resulting model must only have binary variables, have no constraints, and be represented by a polynomial whose degree is at most 2.

We also consider an intermediate representation of the problem based on pseudo-Boolean functions~\cite{boros:2002}, that is,
real polynomials on binary variables of the form
\begin{equation}
    f(\vec{x}) = \sum_{\omega \in \mathcal{P}([n])} c_{\omega} \prod_{j \in \omega} x_{j}
    \label{eq:pbf}
\end{equation}
where $\mathcal{P}([n])$ is the power set of the set $[n] = \{1,\dots,n\}$ and \( f : \set{0, 1}^{n} \to \mathbb{R} \) with \( c_{\omega} \in \mathbb{R} \).
Let \( \mathscr{F} \) be the family of all pseudo-Boolean functions and \( \mathscr{F}^{k} \) its subset containing only polynomials of degree $k$ or less, i.e., \(\mathscr{F}^{k} = \set{f : f \in \mathscr{F}, \deg f \le k} \).
Working with these mathematical objects is a natural choice because, apart from a constant term \(c_{\set{}}\), optimization of \( \mathscr{F}^{2} \) over binary variables is equivalent to QUBO.

Moreover, as previously mentioned at the end of Subsection \ref{sec:miqp}, it is worth reinstating that the final QUBO formulation can be perceived as an approximation to the original model.
However, depending on the reformulation procedures used, the resulting QUBO instance can be a poor representation of the original model.

\subsection*{Variable Encoding}
% When attempting to run constrained, non-binary programs on {QUBO}-powered machines, one must find

When attempting to run constrained or non-binary programs on QUBO-amenable solvers, one must be able to encode the problem's variables and functions as a QUBO instance, which is described by binary variables. 
The conversion from non-binary to binary variables is the first step in the conversion process.
In some cases, the reformulation is exact, for example, when converting bounded integer variables to binary.
However, real variables must go through discretization prior to the aforementioned transformation.

Currently, some well-known encoding techniques are being used in the development of QUBO reformulations.
The most prominent are Binary~\cite{Tamura_Shirai_Katsura_Tanaka_Togawa_2021}, Unary~\cite{Tamura_Shirai_Katsura_Tanaka_Togawa_2021}, One-Hot~\cite{Tamura_Shirai_Katsura_Tanaka_Togawa_2021}, Domain-wall~\cite{Chancellor_2019} and Bounded-coefficient~\cite{karimi2019practical}.
{Additionally, \texttt{ToQUBO.jl} uses a method that we have labeled as Arithmetic Progression~(AP) encoding, which is presented in detail in Appendix \ref{sec:appendix:ap}.}

Each representation has its characteristics, e.g., the number of binary variables~(bits) required, the number of terms in the polynomial expansion produced, and the maximum absolute value $\Delta$ of its coefficients, which is known to affect the quality of the resulting formulation~\cite{quantum_bridge_analytics}.
An asymptotic comparison is presented in Table \ref{tab:encoding}, suggesting that the AP encoding paradigm mitigates the growth of the expansion coefficients while keeping the demand for new variables sublinear.

{These possible variable encoding formulations are set using \texttt{JuMP}'s attribute system, while the binarization technique is controlled by \texttt{ToQUBO.Attributes.VariableEncodingMethod()}.
The default value for this configuration is \texttt{ToQUBO.Encoding.Binary()}, which employs a conservative approach in terms of the number of binary variables added to the final model.
}

{
\renewcommand{\arraystretch}{0.75} %
\begin{table}[ht]
\centering
\begin{threeparttable}[ht]
    \begin{tabular}{|c|c|c|c|c|}
        \hline
        \multirow{2}{*}{Encoding Method} & \multirow{2}{*}{Binary Variables} & \multicolumn{2}{c|}{\# of terms} & \multirow{2}{*}{$\Delta$} \\
        \cline{3-4}
                          &                    & Linear            & Quadratic    & ~              \\
        \hline   
        Binary~\tnote{1}            & $\bigo{\log n}$    & $\bigo{\log n}$   & -            & $\bigo{n}$ \\
        Unary~\tnote{1}           & $\bigo{n}$         & $\bigo{n}$        & -            & $\bigo{1}$     \\
        One-Hot~\tnote{1}           & $\bigo{n}$         & $\bigo{n}$        & $\bigo{n^2}$ & $\bigo{n}$     \\
        Domain-Wall~\tnote{2}       & $\bigo{n}$         & $\bigo{n}$        & $\bigo{n}$   & $\bigo{n}$     \\
        Bounded-Coefficient~\tnote{3} \      & $\bigo{n}$         & $\bigo{n}$        & -   & $\bigo{1}$     \\
        Arithmetic Progression   & $\bigo{\sqrt{n}}$  & $\bigo{\sqrt{n}}$ & -            & $\bigo{\sqrt{n}}$     \\
        \hline
    \end{tabular}
    \begin{tablenotes}
  \footnotesize
  \item[1] \cite{Tamura_Shirai_Katsura_Tanaka_Togawa_2021};
  \item[2] \cite{Chancellor_2019};%
  \item[3] \cite{karimi2019practical};%
\end{tablenotes}
    
\end{threeparttable}
\caption{Comparison between methods for encoding values for \(x \in \{0, \dots, n\}\) into binary variables in \texttt{ToQUBO.jl}}
 \label{tab:encoding}
\end{table}}

\subsection*{Penalty Mechanism}

Representing constraints in an unconstrained problem involves adding penalty terms to the objective function of the final problem.
Constraints can be expressed as functions belonging to feasibility sets, e.g., \( g_{i}(\vec{x}) \leq 0 \iff g_{i}(\vec{x}) \in S_{i}, S_{i} = ( -\infty, 0 ] \).
Every constraint \( g_{i}(\vec{x}) \in S_{i} \) is translated into a penalty function \( \norm{g_{i}(\vec{x})}_{S_{i}} \) with its corresponding penalty factor \( \rho_{i} \).
Under minimization, for example, the accurate portrayal of the feasible set will require each penalty function to be positive if its constraint is violated and zero otherwise.

\begin{equation}
    \begin{array}{>{\displaystyle}r>{\displaystyle}l>{\displaystyle}l}
       \min_{\vec{x}} & f(\vec{x}) \\
        \textrm{s.t.} & g_{i}(\vec{x}) \in S_{i} & \forall i
    \end{array}
    ~~
    \xmapsto{\textrm{Penalization}}
    ~~
    \begin{array}{>{\displaystyle}r>{\displaystyle}l>{\displaystyle}l}
       \min_{\vec{x}} & f(\vec{x}) + \sum_{i} \rho_{i} \norm{g_{i}(\vec{x})}_{S_{i}}
    \end{array}
\end{equation}

Each penalty factor \( \rho_{i} \) should be large enough to ensure that the constraints are satisfied under overall optimality~\cite{quantum_bridge_analytics,Lucas_2014}.
However, if the chosen coefficients are too large, they might have a negative impact on the conditioning of the resulting expression.

\texttt{ToQUBO.jl} implements the mapping and penalization of several constraint types, which is depicted in Section \ref{sec:comparison}.
For example, a common method to embed linear equality constraints will introduce a quadratic expression into the objective function as \( A\,\vec{x} = \vec{b} \) becomes \( (A\,\vec{x} - \vec{b})' (A\,\vec{x} - \vec{b}) \).
However, when dealing with other families of constraints, higher-order penalty functions might arise, and an additional degree-reduction step known as \textit{quadratization} is required for them to fit the QUBO formalism.

\subsection*{Quadratization}

A \textit{quadratization} is a mapping \( \mathcal{Q}: \mathscr{F} \to \mathscr{F}^{2} \) such that
\begin{equation}
    \forall f \in \mathscr{F}, \forall \vec{x} \in \set{0, 1}^{n}, \min_{\vec{w}} \quadratize{f}(\vec{x}; \vec{w}) = f(\vec{x})
\end{equation}
where \( \vec{w} \in \set{0,1}^{m} \) is a vector of auxiliary decision variables.

There are many possible quadratization methods for writing QUBO models from higher-degree pseudo-Boolean functions~\cite{dattani:2019}.
By leveraging \texttt{Julia}'s multiple dispatch paradigm, \texttt{ToQUBO.jl} allows its users to extend the quadratization interface by implementing their degree reduction algorithms. 
For more details on various quadratization schemes, refer to~\cite{dattani:2019}
Currently, two single-term quadratization techniques have already been implemented: one for negative terms and the other for positive terms.

The standard reduction method for negative terms, labeled NTR-KZFD, was introduced by Kolmogorov and Zabih~\cite{kolmogorov2004energy}, and later by Freedman and Drineas~\cite{freedman2005energy}, reduces a single term introducing a single auxiliary decision variable $w$, as follows.
\begin{equation}
    -x_1 x_2 \cdots x_k \xmapsto{\textrm{NTR-KZFD}} (k-1) w - \sum_i x_i w
\end{equation}

Moreover, the default positive term quadratization technique, named PTR-BG, was developed by Boros and Grubner~\cite{boros2014quadratization}.
It is also focused on single terms with a degree higher than two and works with $k-2$ extra decision variables $w_{i}$, where $k$ is the number of variables in the original high-order term, as presented here
\begin{equation}
    x_1 x_2 \cdots x_k \xmapsto{\textrm{PTR-BG}} \left[{
        \sum^{k-2}_{i=1} w_{i} \left({k - i - 1 + x_i + \sum^k_{j=i+1} x_j}\right)
    }\right] + x_{k-1}x_k .
\end{equation}

% At the quadratization step, \texttt{ToQUBO.jl} checks for each term in the objective function whether it has a positive or negative coefficient, as each reduction technique has its dispatch. 
% Then, it filters out quadratic or linear terms, as trying to reduce them would lead to unnecessary computational time consumption.
% After these steps, \texttt{ToQUBO.jl} is ready to apply the required method.
\subsection*{Reformulation Layer}\label{sec:virtual}

\texttt{ToQUBO.jl} implements a mathematical program reformulation layer to provide a transparent interface for the user to model, optimize, and collect results from QUBO sampling runs.
The source program is written as a regular \texttt{JuMP} model, so all the data is cached in the structures defined in \href{https://github.com/jump-dev/MathOptInterface.jl}{\texttt{MathOptInterface}}~(\cite{legat:2021}) to represent variables, constraints, objectives, and additional attributes.
Starting from the well-defined mathematical program formulation stored in the MathOptInterface (\texttt{MOI}) structure, \texttt{ToQUBO.jl} applies the procedures described in Section \ref{sec:compilation} to convert the original problem into a QUBO form that is cached in a new \texttt{MOI} model.
Since the new model is a QUBO defined in \texttt{MOI}, it can be directly sent to any solver that implements the interface and supports the required features, i.e., binary variables and quadratic objectives.
Therefore, we can solve the problem with QUBO sampling machines or even a regular Mixed Integer Quadratic Programming (MIQP) solver capable of handling non-convex objectives, as long as they have an \texttt{MOI} wrapper.
In all these conversions and forwarding steps, the relationship maps are kept in memory so that the user can query the results of the optimization. This process is depicted in Figure \ref{fig:virtual-pipeline-simple}.

The ability to create such layers is an outstanding feature of \texttt{JuMP} and \texttt{MOI} and was a key motivation for selecting this ecosystem to develop a QUBO reformulator.
Previous examples of this strategy are \href{https://github.com/jump-dev/Dualization.jl}{\texttt{Dualization.jl}}~(\cite{bodin:2021}), which receives a \texttt{JuMP} model and provides to the solver its dual, \href{https://github.com/joaquimg/QuadraticToBinary.jl}{\texttt{QuadraticToBinary.jl}}~(\cite{garcia:2021}), which converts quadratic constraints into linear constraints with binary variables, and
% \texttt{DiffOpt.jl}~(\cite{sharma:2022}), which adds the possibility of computing derivatives of solutions of optimization problems with respect to their input data.
\texttt{DisjunctiveProgramming.jl}~(\cite{perez2023disjunctiveprogramming}), which allows the formulation of disjunctive programs and then automatically converts them into MIPs or MINLPs.

\section{Solver Interface}\label{sec:solver}
\texttt{QUBODrivers.jl} provides a common interface for developing bindings that bring QUBO annealing and sampling platforms to the \texttt{JuMP} environment.
By sub-typing MathOptInterface's \texttt{AbstractOptimizer}, the \texttt{AbstractSampler} standard grants the user many relevant features, e.g., querying and sorting multiple results, simple QUBO model format validation and adjustable attributes for fine-tuning runtime execution.
By employing \texttt{QUBODrivers.jl}, the user can define a new MOI-compliant solver interface, thus providing access to different hardware architectures and sampling algorithms with ease.

As presented in Listing \ref{code:drivers}, everything begins with the \texttt{@setup} macro, whose body contains all the necessary settings to specify a new optimizer. 
When declaring such properties, one should first fill in the solver's \texttt{name}, \texttt{sense}, \texttt{domain}, and \texttt{version}, followed by the solver-specific parameters, which belong to the \texttt{attribute} block.
This macro leverages \texttt{Julia}'s meta-programming capabilities to circumvent a considerable amount of repetitive code and isolate the user from most of MOI's internal details.

By overloading the \texttt{sample} function, the QUBO model parsed from \texttt{JuMP} is then used to sample and communicate back solutions from an arbitrary underlying procedure.
The return of this function must be a \texttt{SampleSet}, a special collection designed to provide fast queries and serve as input to analytical tools.

\begin{jllisting}[caption = {Example of \texttt{QUBODrivers.jl} usage to wrap QUBO solvers}, label = code:drivers, captionpos=b]
module RandomSampler

import Random
import MathOptInterface as MOI
import QUBODrivers, QUBOTools

QUBODrivers.@setup Optimizer begin
    name = "Random Sampler"
    sense = :min
    domain = :bool
    version = v"1.0.0"
    attributes = begin
        RandomSeed["seed"]::Union{Integer,Nothing} = nothing
        NumberOfReads["num_reads"]::Integer        = 1_000
    end
end

function QUBODrivers.sample(sampler::Optimizer{T}) where {T}
    # Retrieve Raw Model
    Q, a, b = QUBODrivers.qubo(sampler, Dict)

    # Retrieve User-Defined Attributes
    n         = MOI.get(sampler, MOI.NumberOfVariables())
    num_reads = MOI.get(sampler, RandomSampler.NumberOfReads())
    seed      = MOI.get(sampler, RandomSampler.RandomSeed())
    Random.seed!(seed)

    # Sample Random States
    samples = QUBOTools.Sample{T,Int}[]
    
    for _ in 1:num_reads
        x = Random.rand((0, 1), n)
        y = QUBOTools.value(Q, x, a, b)
        push!(samples, QUBOTools.Sample{T,Int}(x, y))
    end

    # Return Solution
    return QUBOTools.SampleSet{T,Int}(samples)
end
end # module
\end{jllisting}

\section{Tools for QUBO} \label{sec:format-tools}

Last but not least, we present \texttt{QUBOTools.jl}, a library containing interface definitions and core functionality for managing and analyzing both QUBO models and their solutions.
Its development is focused on providing
\begin{itemize}
    \item Fast and reliable I/O that allows conversion between well-known file formats for QUBO and Ising models
    \item Generic modeling back-end utilities for powering other applications, including specialized reference implementations for conceptual data structures
    \item Analytical framework equipped with data queries, conditioning, performance evaluation metrics, and plot recipes
    \item Routines for random instance generation (see Section \ref{sec:next-steps} - Next Steps)
\end{itemize} 

\subsection{File Formats}

As discussed above, Quantum Operations Research has many agents, including hardware manufacturers, software solution providers, and mathematical optimization experts.
An ecosystem as diverse and dynamic as this also has the drawback of having many different file formats and interfaces for specifying and running QUBO models.
To address this issue, \texttt{QUBOTools.jl} provides an extensible I/O library with built-in \textit{codecs}, making it a file converter compatible with the most widely used formats.

The \textit{bqpjson} variant, for instance, is a platform-agnostic \texttt{JSON} specification developed by Los Alamos National Laboratory's Advanced Network Science Initiative to incorporate additional aspects of the problem, including related solutions and metadata from the target hardware~\cite{bqjson:2020}.
Another example, the \textit{QUBO} format, is defined according to D-Wave's applications for quantum and simulated annealing, but is also configurable to establish communication with MQLib.
The \textit{Qubist} type, also standardized by D-Wave, is one of the simplest file standards available, with no support for metadata entries or any other sophisticated feature.

However, once able to consistently interpret \textit{Qubist} and a few other file formats, it is possible to gain instant access to many problem databases produced by a wide range of groups~\cite{hen:2019,kowalsky:2022}.
Conversely, being capable of writing in several QUBO dialects allows one to leverage the most prominent solvers available when building new applications.
By operating in both ways, \texttt{QUBOTools.jl} establishes the technical foundation for performing intricate tasks such as benchmarking heterogeneous hardware and formulation analysis. 
{
The forementioned file formats are documented with examples in \texttt{QUBOTools.jl} documentation~\footnote{\texttt{QUBOTools.jl} documentation for file formats: \href{https://juliaqubo.github.io/QUBOTools.jl/dev/manual/5-formats/}{juliaqubo.github.io/QUBOTools.jl/dev/manual/5-formats/}}}

\subsection{Data Structures}

A project like \texttt{QUBO.jl} spans multiple activities within what could be called the "application layer" of Quantum Optimization Systems.
It plays the role of a user-friendly platform by exposing high-level methods and also internally operating closer to the bare metal, e.g., accounting for hardware specifics when needed.
Taking into account the different descriptions of the models and solver architectures, \texttt{QUBOTools.jl} defines a blueprint to coordinate the integration of \texttt{QUBO.jl}'s components.

Harnessing \texttt{Julia}'s powerful type system, a manifold of generic methods is delineated over a few conceptual data structures.
This abstract interface puts the ecosystem's overall design into perspective, focusing on a functional configuration rather than a static set of tools.
With that said, \texttt{QUBOTools.jl} delivers a protocol to work with models, solvers, and solutions.

Moreover, each envisioned data structure has a self-contained, ready-to-use reference implementation.
The main goal of this approach is to offer the building blocks for other applications, including \texttt{QUBO.jl} itself.
By exploiting both abstract and concrete functionalities, one might be able to specialize portions of software behavior, either by extending the proposed interface or by the composition of its main artifacts.

\subsection{Analytical Tools}

Given the variety of quantum hardware in the NISQ Era, it is still unknown which candidates will prevail as the field moves towards Fault-Tolerant QC.
Furthermore, each device can perform better with a specific QUBO formulation, taking into account variable encoding, constraint mapping, and quadratization techniques. 
Under these circumstances, benchmarking NISQ systems and QUBO models in the context of Quantum Optimization provides a guiding framework to assess the current state of these emerging solvers.

Analytical tools were bundled into the package's library to address this concern.
These tools provide access to the conditioning and structure of the model and the quality of the sampled solutions.
When making comparisons between solvers, \texttt{QUBOTools.jl} can gauge time-related indicators, inspect the success rate of a set of runs, and evaluate composite metrics such as the time-to-solution~(TTS)~\cite{king2015benchmarking}.
Also, in light of model conditioning, one can withdraw statistics from the formulation amongst model density and the connectivity arising from its quadratic relations.
The latter could also be helpful in choosing between reformulation techniques in the context of \texttt{ToQUBO.jl}.

\subsubsection*{Visualization}

Adding to the data analysis toolkit, \texttt{QUBOTools.jl} also provides a set of predefined visualization schemes.
It leverages \texttt{Plots.jl}'s~\cite{plots_jl} submodule \texttt{RecipesBase.jl} to define custom plot instructions, considering the different perspectives through which a model or solver can be visually examined.
These recipes are a series of drawing statements and are not tied to any specific plotting library. Thus, they considerably reduce the overhead by delaying or avoiding the installation of extra dependencies.

\begin{figure}[!ht]
\centering
\subfloat[]{%
\resizebox*{8cm}{!}{\includegraphics[align=c]{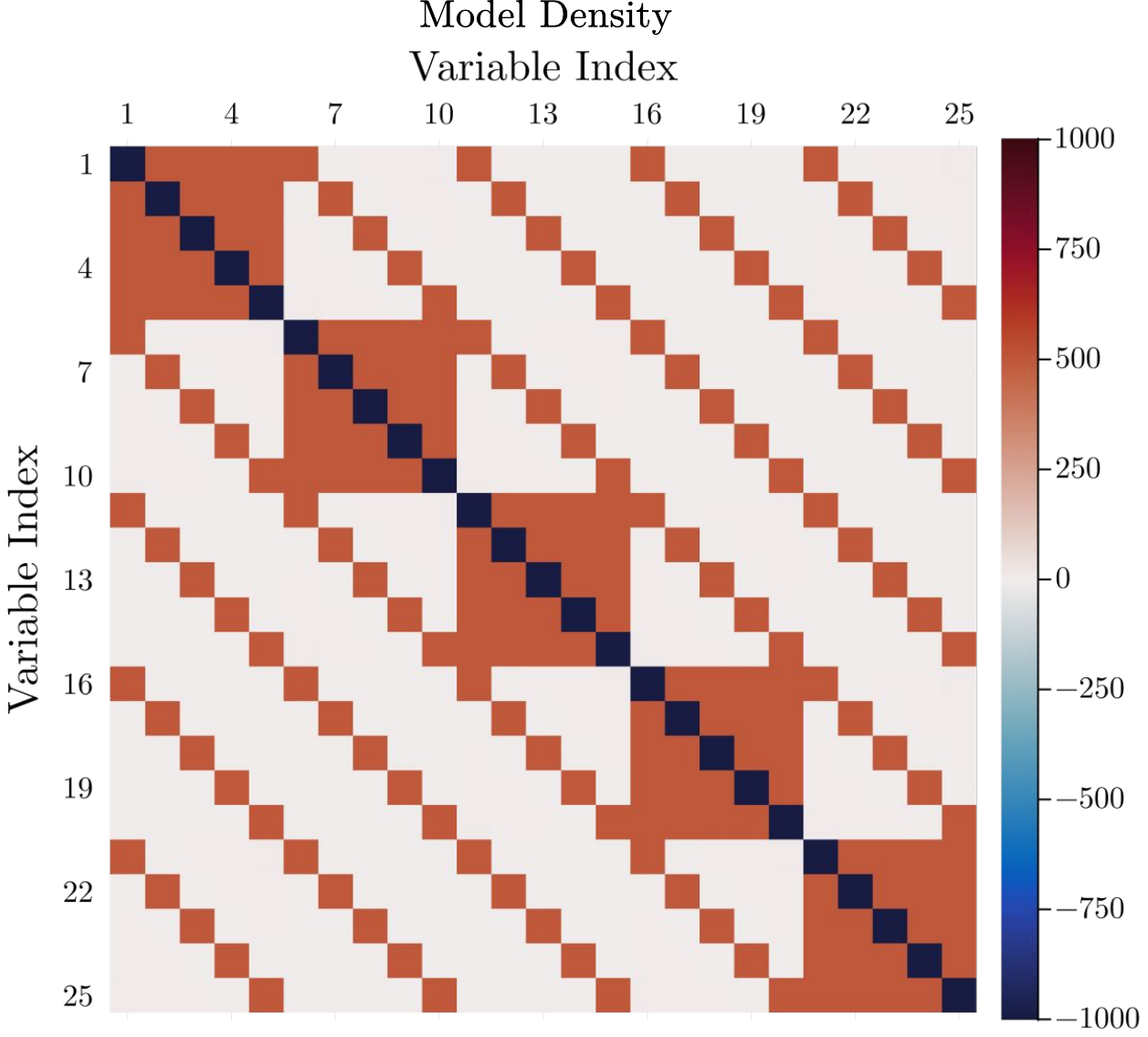}}}\\
\subfloat[]{%
\resizebox*{9cm}{!}{\includegraphics[align=c]{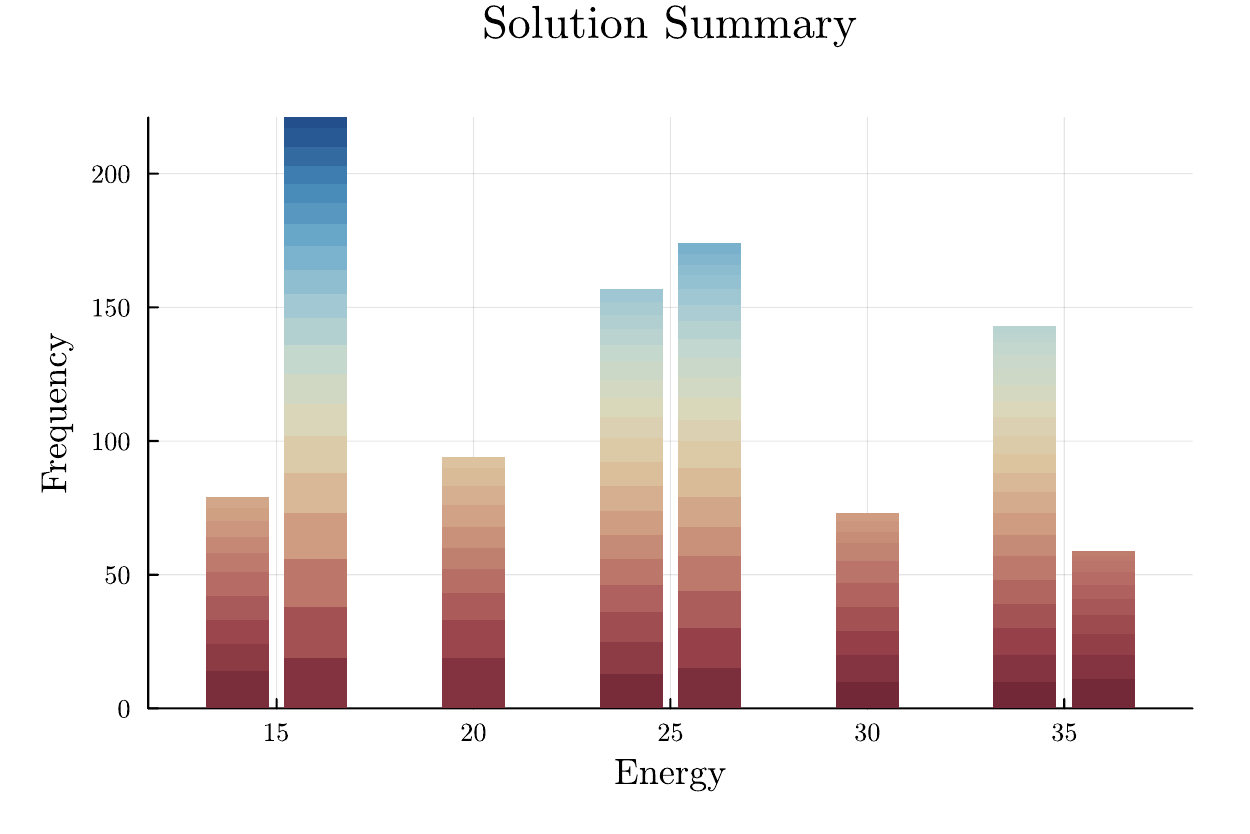}}}
\caption{\texttt{QUBOTools.jl}'s plot recipes in action: model density~(Subfigure a), from a TSP instance with 5 cities~(25 variables) and sampling frequency by energy value~(Subfigure b), using the \texttt{PySA} solver, {where samples with the same variable values are grouped by color within the same column (the same color in different columns describes no relation).} \\}
    \label{fig:plots}
\end{figure}

By design, \texttt{QUBOTools.jl} aims to provide the most widely used figures for quantum optimization analysis, extending \texttt{Plots.jl} built-in methods.
Currently, plotting recipes for examining the outcomes' sampling distribution and the model's matrix density have already been implemented.
On its left side, Figure \ref{fig:plots} depicts an example of a Traveling Salesperson Problem (TSP), where the color saturation is proportional to the magnitude of each coefficient.
On the right side, samples from running the TSP model are sorted by their respective objective values, and the number of reads from each configuration is given by their heights.
The results were gathered using the NASA Parallel Tempering open-source solver \texttt{PySA}~\cite{mandra2018deceptive, Mandra_PySA_Fast_Simulated_2023}, and samples from different states but with similar energy values were stacked together.

% \begin{figure}[h]
% \centering
% \subfloat[]{%
% \resizebox*{5cm}{!}{\includegraphics[align=c]{figures/model_density_tsp.pdf}}}\hspace{1pt}
% \subfloat[]{%
% \resizebox*{8cm}{!}{\includegraphics[align=c]{figures/energy_frequency_tsp.pdf}}}
% \caption{\texttt{QUBOTools.jl}'s plot recipes in action: model density~(left), from a TSP instance with 5 cities~(25 variables) and sampling frequency by energy value~(right), using the \texttt{PySA} solver, {where samples with the same variable values are grouped by color within the same column (the same color in different columns describes no relation).} \\}
%     \label{fig:plots}
% \end{figure}

Both conditioning metrics and visualization recipes were built on the basis of the abstract data structure layout mentioned above.
A consistent set of methods for querying models and solutions is one of the fundamental building blocks of the entire \texttt{QUBO.jl} software stack.
By wrapping \texttt{JuMP} models and their solution sets with constructs that adhere to its interface, \texttt{QUBO.jl} allows users to leverage the analytical toolbox presented seamlessly.
\section{The State of Software for Quantum Computing}\label{sec:julia_env}

Evaluating the current scenario for software packages for Quantum and Quantum-Inspired Optimization with QUBO applications, it is noticeable that \texttt{Python} is the most popular language for their development.
Some popular examples are \href{https://amplify.fixstars.com/en/}{\texttt{Amplify}}~\cite{fixstars}, \href{https://github.com/recruit-communications/pyqubo}{\texttt{PyQUBO}}~\cite{pyqubo:2021}, \href{https://github.com/jtiosue/qubovert}{\texttt{qubovert}}~\cite{qubovert:2022} and DWave's \href{https://github.com/dwavesystems/dimod}{\texttt{dimod}}~(a module of the Ocean SDK package).
In addition, IBM's \texttt{Qiskit}~\cite{Qiskit}, Xanadu's \href{https://pennylane.ai/}{\texttt{PennyLane}}~\cite{pennylane} and EntropicaLab's \href{https://github.com/entropicalabs/openqaoa}{\texttt{OpenQAOA}}~\cite{openqaoa} also stand out.
Moreover, \texttt{SATyrus}~\cite{satyrus:2007} is a logic-oriented modeling platform capable of generating QUBO models from descriptions of satisfiability (SAT)-like problems.
Its latest implementation, \texttt{SATyrus III}~\cite{satyrus:2022}, served as an initial inspiration for some of the design choices behind \texttt{QUBO.jl}.

Even in light of this panorama, the \texttt{Julia} ecosystem for QC has been gaining momentum, with Amazon recently releasing a \href{https://github.com/awslabs/Braket.jl}{\texttt{Braket.jl}}~\cite{braketjl} version that supports communication with its QC environment, AWS Braket. 
In addition, QuEra, a quantum hardware company, which is also available on AWS Braket, has published \href{https://github.com/QuEraComputing/Bloqade.jl}{\texttt{Bloqade.jl}}~\cite{bloqade2023quera}, a software to interface with its systems in \texttt{Julia}. 
{  Los Alamos Advanced Network Science Initiative has also published the \href{https://github.com/lanl-ansi/QuantumAnnealing.jl}{\texttt{QuantumAnnealing.jl}}~\cite{quantumannealing:2022, quantum_annealing_jl_paper} package for simulation of quantum annealing algorithms.}
Other examples include the \href{https://github.com/QuantumBFS/Yao.jl}{\texttt{Yao.jl}}~\cite{YaoFramework2019} framework for experimenting with quantum algorithms, \href{https://github.com/harshangrjn/QuantumCircuitOpt.jl}{\texttt{QuantumCircuitOpt.jl}}~\cite{QCOpt_SC2021} for optimizing quantum circuit design.
\texttt{ToQ.jl}~\cite{toq_jl} is a library provided to interface with D-Wave's quantum annealers, later continued as \href{https://github.com/omalled/ThreeQ.jl}{\texttt{ThreeQ.jl}}.
{Additionally, the packages \texttt{QAOA.jl}~\cite{qaoa_jl} and \texttt{JuliQAOA.jl}~\cite{juliqaoa} present implementations of the QAOA algorithm.}

As for Quantum Optimization,
\texttt{Julia} provides essential features on top of the above-mentioned growing environment, as it was designed for high-performance computing \cite{julia:2017}.
Hence, it allows for the development of efficient algorithms and avoids the two-language problem that occurs when part of an interpreted language, e.g., \texttt{Python}, code must be rewritten in a compiled one, e.g., \texttt{C++}, when aiming for performance.
This has been the case for \href{https://github.com/recruit-communications/pyqubo}{\texttt{PyQUBO}}~\cite{pyqubo:2021} and \href{https://amplify.fixstars.com/en/}{\texttt{Amplify}}~\cite{fixstars}, where the core of the reformulator to QUBO problems had to be rewritten from \texttt{Python} to \texttt{C++}.
Also, \texttt{JuMP} is an open-source library used by the Operations Research and Optimization communities because it supports a broad range of Optimization classes, is efficient, extensible, and user-friendly.

Moreover, programming in \texttt{Julia} does not exclude one from using quantum computers, e.g.,  using \texttt{Qiskit} to access IBM's hardware.
Building Python package wrappers with tools such as \href{https://github.com/cjdoris/Pythoncall.jl}{\texttt{PythonCall}}~\cite{PythonCall} allows exploring the benefits of \texttt{Julia} and the wide variety of quantum software packages in \texttt{Python}.
This review of the software status for quantum computing motivated us to write our tools in \texttt{Julia}, leveraging the strengths of the programming language and the existing package environment for operations research and quantum computing.

\section{Comparison between existing tools}\label{sec:comparison}

As mentioned in the previous section, there are a few available tools for working with QUBO formulations, and the most relevant lie within the Python environment, such as \texttt{Amplify}, \texttt{PyQUBO}, \texttt{qubovert}, \texttt{dimod}, \texttt{Qiskit} and \texttt{OpenQAOA}. 
These modules differ mainly in their capabilities, such as Automatic Variable Encoding and Constraint Mapping coverage. From this perspective, we have listed these features for each tool in Table \ref{tab:tools}, showing that \texttt{ToQUBO.jl} has all types of encodings and constraints covered by the mentioned packages. 

{
    \renewcommand{\arraystretch}{1.0} %
    \begin{table}[ht]
    \small
    \centering
    \def\YES{$\blacksquare$}
    \def\YESBUT#1{\phantom{\tnote{#1}}~\YES~\tnote{#1}}
    \def\NOBUT#1{}%{\phantom{\tnote{#1}}~\phantom{\YES}~\tnote{#1}}
    \begin{threeparttable}
   
    \begin{tabularx}{\linewidth}{|C|>{\leavevmode\color{PSRGOLD}\bfseries}c|>{\leavevmode\color{PSRBLUE}\bfseries}c|>{\leavevmode\color{PSRBLUE}\bfseries}c|>{\leavevmode\color{PSRBLUE}\bfseries}c|>{\leavevmode\color{PSRBLUE}\bfseries}c|>{\leavevmode\color{PSRBLUE}\bfseries}c|>{\leavevmode\color{PSRBLUE}\bfseries}c|}
        \hline
   & {\footnotesize\texttt{ToQUBO.jl}} & {\footnotesize\texttt{PyQUBO}} & {\footnotesize\texttt{qubovert}} & {\footnotesize\texttt{Qiskit}} & {\footnotesize\texttt{OpenQAOA}} & {\footnotesize\texttt{Amplify}} & {\footnotesize\texttt{dimod}}\\

        \hline
        \multicolumn{8}{|c|}{Automatic variable encoding methods implemented}                                                                                                     \\
        \hline
        Binary~\tnote{1}                                     & \YESBUT{$\bigstar$}                 & \YES              & \YES               & ~~                  & ~~ & \YESBUT{!}     & \YES       \\
        Unary~\tnote{1}                                      & \YESBUT{$\bigstar$}                 & \YES              & \YES               & \YES                  & \YES & \YESBUT{!}       & ~~      \\
        One-Hot~\tnote{1}                                    & \YESBUT{$\bigstar$}                 & \YES              & ~~               & ~~                  & ~~ & ~~       & ~~     \\
        Domain-Wall~\tnote{2}                                & \YESBUT{$\bigstar$}                 & \YES              & ~~               & ~~                  & ~~ & ~~     & ~~       \\
        Bounded-Coefficient~\tnote{3}                        & \YESBUT{$\bigstar$}                 & \YES              & ~~               & ~~                  & ~~ & ~~       & ~~     \\
        Arithmetic Progression                                   & \YESBUT{$\bigstar$}                 & ~~              & ~~               & ~~                  & ~~ & \YESBUT{!}     & ~~       \\
        \hline  
        \multicolumn{8}{|c|}{Supported automatic constraint reformulation}                                                                      \\
        \hline  
        \( \vec{a}'\vec{x} \le b \)                        & \YES                 & ~~              & \YES               & \YES                    & \YES & \YES        & \YES   \\
        \( \vec{a}'\vec{x} = b \)                          & \YES                 & ~~              & \YES               & \YES                    & \YES & \YES        & \YES    \\
        \( \vec{x}' \vec{Q}\,\vec{x} + \vec{a}'\vec{x} \le b \)   & \YES                 & ~~    & ~~               & ~~                    & ~~ & \YES          & \YES   \\
        \( \vec{x}' \vec{Q}\,\vec{x} + \vec{a}'\vec{x} = b \)     & \YES                 & ~~    & ~~               & ~~                    & ~~ & \YES        & \YES     \\
        SOS1~\tnote{4}                                              & \YES                 & ~~              & ~~               & ~~                    & ~~ & ~~       & ~~      \\
        \( \bigwedge_{i} x_{i} \)~\tnote{$\dagger$}                                               & \NOBUT{$\blacklozenge$}                & ~~              & \YES               & ~~                    & ~~ & \YES     & ~~        \\
        \( \bigvee_{i} x_{i} \)~\tnote{$\dagger$}                                               & \NOBUT{$\blacklozenge$}                 & ~~              & \YES               & ~~                    & ~~ & \YES     & ~~        \\
        \( \bigoplus_{i} x_{i} \)~\tnote{$\dagger$}                                               & \NOBUT{$\blacklozenge$}                 & ~~              & ~~               & ~~                    & ~~ & \YES       & ~~      \\
        \hline
    \end{tabularx}
    \begin{tablenotes}
      \footnotesize
      \item[1] \cite{Tamura_Shirai_Katsura_Tanaka_Togawa_2021};
      \item[2] \cite{Chancellor_2019};%
      \item[3] \cite{karimi2019practical};%
      \item[4] Special Ordered Set of Type 1~\cite{beale1970special};
      \item[$\dagger$] For logical constraints, it is assumed that \( x_{i} \in \set{0, 1} \). \\
      \item[$\bigstar$] \texttt{ToQUBO.jl} is the only platform in the board whose automatic encoding routines are also applicable to continuous variables. Furthermore, its \textit{Bounded-Coefficient} implementation allows the technique to be used not just with the \textit{Binary} encoding, as proposed in the original paper, but also with the \textit{Unary} and \textit{Arithmetic Progression} modes.\\
      \item[!] \texttt{Amplify} only supports these encoding methods when introducing integer slack variables for mapping inequality constraints. They are not applicable to regular variables, as their models only accept binary or spin sites so far.
    \end{tablenotes}
    \end{threeparttable}
     \caption{List of supported variable encoding and constraint penalization methods across QUBO modeling platforms}
    \label{tab:tools}
    \end{table}
}

As presented earlier, \texttt{QUBODrivers.jl} enables users to interface with QUBO samplers and define new solvers, which work as wrappers, gaining access to different hardware.
As of today, \texttt{QUBODrivers.jl} provides three simple samplers for testing purposes: {\textit{ExactSampler}}, {\textit{RandomSampler}} and {\textit{IdentitySampler}}. 
First, \textit{ExactSampler} performs an exhaustive search over all possible states, while \textit{RandomSampler} assigns a boolean value to each variable, given a probability distribution.
Finally, \textit{IdentitySampler} just selects the state provided as warm-start as the solution.

In addition, we have already released several open-source packager wrappers using \texttt{QUBODrivers.jl}: \href{https://github.com/JuliaQUBO/DWave.jl}{\texttt{DWave.jl}}\cite{dwavejl:2022},
\href{https://github.com/JuliaQUBO/DWaveNeal.jl}{\texttt{DWaveNeal.jl}}\cite{dwavenealjl:2022}, 
\href{https://github.com/JuliaQUBO/MQLib.jl}{\texttt{MQLib.jl}}\cite{mqlib:2022}, 
\href{https://github.com/JuliaQUBO/QuantumAnnealingInterface.jl}{\texttt{QuantumAnnealingInterface.jl}}\cite{quantumannealinginterfacejl:2022},
\href{https://github.com/JuliaQUBO/QiskitOpt.jl}{\texttt{QiskitOpt.jl}}~\cite{qiskitoptjl},
\href{https://github.com/JuliaQUBO/CIMOptimizer.jl}{\texttt{CIMOptimizer.jl}}~\cite{cim-optimizer:2022,cimoptimizerjl:2023},
and \href{https://github.com/JuliaQUBO/PySA.jl}{\texttt{PySA.jl}}~\cite{mandra2018deceptive,pysajl:2023}.
\texttt{DWave.jl} and \texttt{DWaveNeal.jl} were developed as D-Wave's Quantum and Simulated Annealing interfaces, respectively.
\texttt{MQLib.jl} is a wrapper for a heuristics library for QUBO problems~\cite{mqlib:2018} and \texttt{QuantumAnnealingInterface.jl} is an interface for Los Alamos National Laboratory's \texttt{QuantumAnnealing.jl}~\cite{quantumannealing:2022, quantum_annealing_jl_paper}. Finally, \texttt{QiskitOpt.jl} allows users to send their \texttt{JuMP} instances to run VQE and QAOA algorithms.

On the other hand, \texttt{PyQUBO}, \texttt{qubovert} and \texttt{dimod} only offer access to D-Wave's Simulated and Quantum Annealing hardware as external samplers and do not allow users to define new solvers.
Furthermore, in addition to having access to D-Wave's systems, \texttt{Amplify} provides an interface with Fujitsu's Digital Annealer, Toshiba's Bifurcation Machine, Hitachi's CMOS Annealing Machine, IBM's gate-based computers, and other Simulated Annealing machines.
Finally, \texttt{Qiskit} and \texttt{OpenQAOA} allow users to run QUBO formulations on QAOA, VQE, and some other Variational Quantum algorithms that can be executed in IBM's hardware.

It is noticeable that \texttt{Amplify} and \texttt{QUBO.jl} interface with a similar number of solvers hardware-wise, as mentioned in Section \ref{sec:solver}.
However, by enhancing Julia's composability and interoperability features, \texttt{QUBODrivers.jl} users can easily define new solver wrappers. 
As the \texttt{QUBO.jl} ecosystem gains momentum, we expect that our sampler coverage will increase, as it is an open-source and extensible environment designed for algorithm researchers, hardware manufacturers, and enthusiasts in general.
This comparison is presented in Table \ref{tab:interfaces}.

{
    \renewcommand{\arraystretch}{1.5} %
    \begin{table}[ht]
    \scriptsize
    \centering
    \def\S{$\blacksquare$}
    \def\T#1{\phantom{\tnote{#1}}~\S~\tnote{#1}}
    \begin{threeparttable}
    \begin{tabularx}{\linewidth}{|c|C|>{\leavevmode\color{PSRGOLD}\bfseries}c|>{\leavevmode\color{PSRBLUE}\bfseries}c|>{\leavevmode\color{PSRBLUE}\bfseries}c|>{\leavevmode\color{PSRBLUE}\bfseries}c|>{\leavevmode\color{PSRBLUE}\bfseries}c|>{\leavevmode\color{PSRBLUE}\bfseries}c|>{\leavevmode\color{PSRBLUE}\bfseries}c|}
        \hline
        Sampler & Hardware &
\rotatebox{60}{\texttt{QUBODrivers.jl}} &
\rotatebox{60}{\texttt{PyQUBO}} &
\rotatebox{60}{\texttt{qubovert}} &
\rotatebox{60}{\texttt{Qiskit}} &
\rotatebox{60}{\texttt{OpenQAOA}} &
\rotatebox{60}{\texttt{Amplify}} &
\rotatebox{60}{\texttt{dimod}} \\
        \hline  
        Simulated Annealing                   & CPU/ GPU      & \T{3}     & \S         & \S            & ~~        & ~~ & \S  & \S\\ 
        Quantum Annealing                & QA\tnote{1}  & \T{4}     & \S         & ~~                  &  \T{10}       & ~~ & \S  & \S\\ 
        \makecell{Simulated \\ Quantum Annealing}                          & CPU/ GPU      & \T{5}     & ~~         & ~~            & ~~   & ~~  & \S & ~~  \\ 
        VQE                                   & QGC~\tnote{2}         & \T{6}     & ~~         & ~~             & \S       & \S & ~~ & ~~  \\ 
        QAOA                                  & QGC~\tnote{2}          & \T{6}     & ~~         & ~~             & \S       & \S & \S & ~~\\ % 
        \makecell{\texttt{MQLib} \\ Heuristics Library}                   & CPU          & \T{7}     & ~~            & ~~            & ~~   & ~~   & ~~ & ~~   \\  
        Parallel Tempering                & CPU          & \T{8}                 & ~~         & ~~ & ~~ & ~~ & ~~   & ~~    \\
        \makecell{Simulated \\ Bifurcation Machine}       &  GPU/ FPGA     & ~~                 & ~~         & ~~ & ~~ & ~~ & \S  & ~~       \\ 
        \makecell{Simulated Coherent \\ Ising Machine} & CPU/ GPU         & \T{9}                 & ~~         & ~~ & ~~ & ~~ & ~~    & ~~    \\ 
        Digital Annealing         & ASIC         & ~~                 & ~~         & ~~ & ~~ & ~~ & \S  & ~~     \\ 
        CMOS Annealing            & CMOS & ~~                 & ~~         & ~~ & ~~ & ~~ & \S  & ~~     \\
        \hline
    \end{tabularx}
    \begin{tablenotes}
      \footnotesize
      \item[1] {Quantum Annealer};%
      \item[2] {Quantum Gate Circuit};%
      \item[3] {\texttt{DWave.jl}~\cite{dwavejl:2022}};%
      \item[4] {\texttt{DWaveNeal.jl}~\cite{dwavenealjl:2022}};%
      \item[5] {\texttt{QuantumAnnealingInterface.jl}~\cite{quantumannealing:2022,quantumannealinginterfacejl:2022}};%
      \item[6]{\texttt{QiskitOpt.jl}~\cite{qiskitoptjl}};%
      \item[7]{\texttt{MQLib.jl}~\cite{mqlib:2018,mqlib:2022}};%
      \item[8]{\texttt{PySA.jl}~\cite{pysajl:2023}};%
      \item[9]{\texttt{CIMOptimizer.jl}~\cite{cim-optimizer:2022,cimoptimizerjl:2023}};%
      \item[10]{\texttt{Dwave-Qiskit-plugin}~\cite{dwaveqiskitplugin}};%
    \end{tablenotes}
    \end{threeparttable}
    \caption{External QUBO-amenable samplers coverage for each package}\label{tab:interfaces}
    \end{table}
}

\begin{figure}[!ht]
    \centering
    \includegraphics[width=\linewidth]{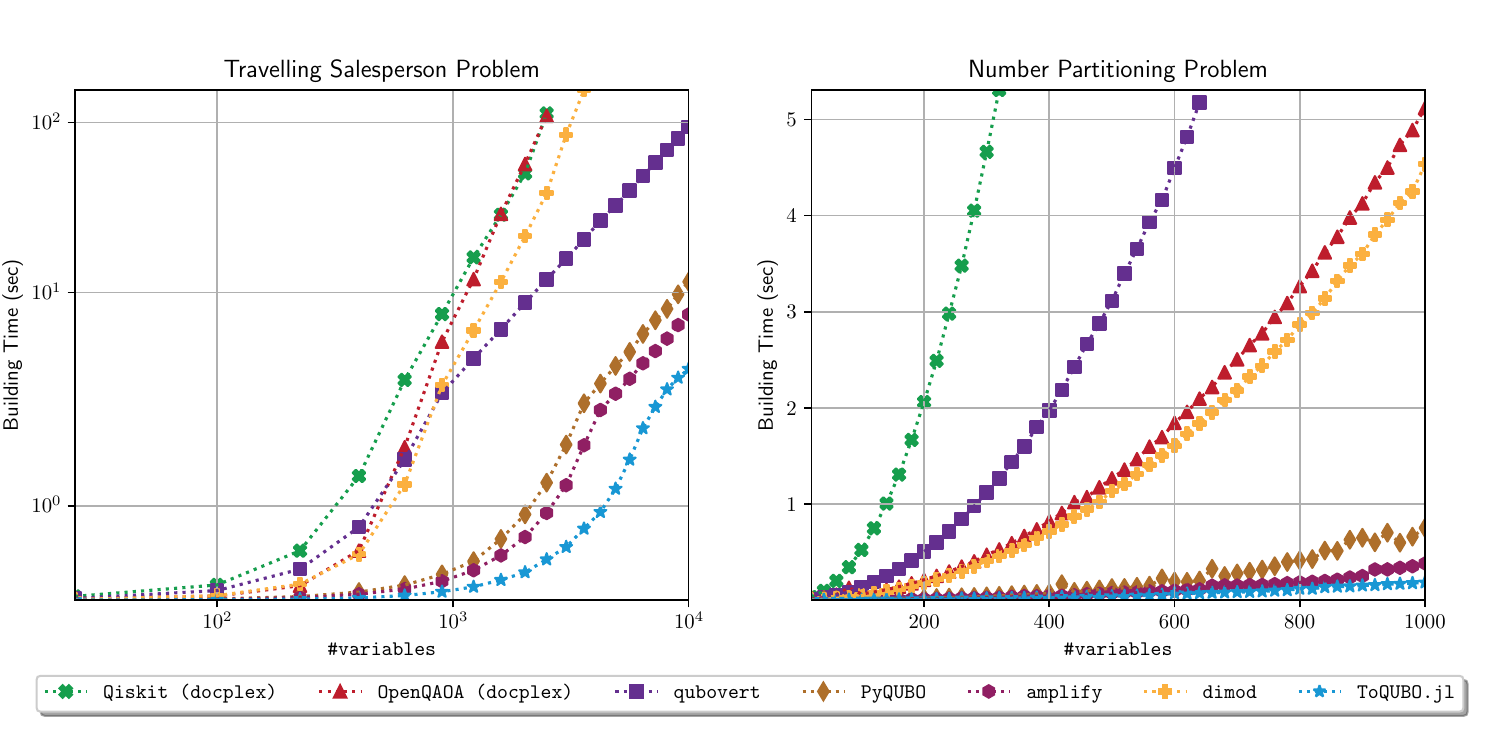}
    \caption{%
    Benchmark results for the Travelling Salesperson Problem~(TSP) and Number Partitioning Problem~(NPP) QUBO building time. Each TSP instance is comprised of \( \sqrt{\text{\texttt{\#variables}}} \) cities. The partitioned sets in the NPP are of the form \( [\text{\texttt{\#variables}}] = \left\lbrace{1, \dots, \text{\texttt{\#variables}}}\right\rbrace \).
    }
    \label{fig:benchmark}
\end{figure}

Moreover, another method to compare QUBO formulation packages is to analyze their execution time to build an actual QUBO model. 
Although at first it is expected for a Julia package to have a better performance than a Python one, considering their design differences, such as the fact that Julia is a compiled language and Python is an interpreted one, \texttt{PyQUBO}'s and \texttt{Amplify}'s backends are written in C++/C.
We have decided to follow \texttt{PyQUBO}'s benchmarking framework, presented in its documentation, where many instances of the Traveling Salesperson Problem (TSP) and the Number Partitioning Problem (NPP) are formulated with a varying number of variables~\cite{Lucas_2014}.
As presented in Figure \ref{fig:benchmark}, we evaluate the time to construct a QUBO model from 25 to 10000 variables (5 to 100 cities) for the TSP and from 5 to 1000 variables for the NPP. It is important to mention that for the \texttt{Qiskit} and \texttt{OpenQAOA} iterations, we create a QUBO from a \texttt{docplex} model~\cite{docplex:2015}.

Examining the results for building a QUBO model from TSP and NPP instances~(see {Figure \ref{fig:benchmark}}) together with the data from Table \ref{tab:tools}, it is possible to conclude that \texttt{ToQUBO.jl} stands out from its main competitors.
In contrast with the other packages, \texttt{ToQUBO.jl} contains more features, providing distinguished compatibility with some problems.
Besides, \texttt{ToQUBO.jl}'s execution time to build a QUBO model is very competitive, being faster than \texttt{PyQUBO} and greatly faster than \texttt{qubovert}, \texttt{dimod}, \texttt{Qiskit}, and \texttt{OpenQAOA}.

% Furthermore, in the perspective of programming experience, being developed as an \texttt{MOI} layer, \texttt{ToQUBO.jl} allows users to use \texttt{JuMP} for modeling, which is a well documented and efficient tool for a broad range of Optimization problem classes.
% On the other hand, \texttt{qubovert} and \texttt{PyQUBO} define their own taylor-made syntax with the sole objective of optimizing QUBO problems.
% We illustrate this advantage with Listings \ref{code:toqubo}, \ref{code:pyqubo}, and \ref{code:qubovert}, used to perform TSP benchmarking.
% Note that the only non-standard \texttt{JuMP} function used is \texttt{ToQUBO.qubo(model)}, as \texttt{ToQUBO.Optimizer} can be regarded as a \texttt{JuMP} solver.

% \input{code/toqubo.tex}
% \input{code/pyqubo.tex}
% \input{code/qubovert.tex}

\texttt{ToQUBO.jl}'s performance and features, united with \texttt{JuMP}'s mathematical programming ecosystem and \texttt{QUBODrivers.jl}, allow \texttt{Julia} users to experiment using quantum and classical hardware-accelerated QUBO solvers in a seamless experience, just as they would do with a standard \texttt{JuMP} solver.

\section{Conclusion}\label{sec:conclusions}

\texttt{QUBO.jl} is presented as an extensive tool for working with QUBO formulations and interfacing with their solvers, standing out from other packages in the Quantum Optimization field.
Taking advantage of \texttt{ToQUBO.jl}'s fast reformulation routines, users can submit \texttt{JuMP} models to emergent QUBO devices as they would do with any other solver.
Moreover, leveraging different architectures, each depending on specific communication protocols, is simplified with \texttt{QUBOTools.jl}.
In addition, since all three modules are open source, anyone can contribute to their development and even define new solvers with \texttt{QUBODrivers.jl}, as discussed in Section \ref{sec:solver}, or create new file conversions with \texttt{QUBOTools.jl}.

Providing all three packages in \texttt{QUBO.jl} as a bundle makes them easy to integrate into optimization specialists' work, as it would be just an abstraction layer between their \texttt{JuMP} instances and quantum or classical sampling devices. 
As our project develops, both academia and industry will be able to use \texttt{QUBO.jl} as a platform for quantum and hardware-accelerated methods, either by connecting their research pipelines to several available QUBO solvers or by publishing their own experimental hardware and algorithms to this growing ecosystem.
\texttt{QUBO.jl} addresses the gap between classical optimization and hardware-accelerated QUBO solvers~(e.g., coherent Ising machines, quantum annealers, and variational quantum algorithms), because it has almost no prerequisites for \texttt{JuMP} developers. 
Moreover, users of other algebraic modeling languages for Optimization only need to learn \texttt{JuMP} and not much about these emerging platforms.

\subsection{Next Steps}\label{sec:next-steps}

We plan to expand the constraint mapping capabilities for future versions of our packages and consider new variable encoding techniques while still fulfilling competitive performance goals. 
More specifically, it is expected that, by employing the expertise earned from \texttt{SATyrus}'s development~\cite{satyrus:2022}, it will be possible to enhance \texttt{ToQUBO.jl} with support for logical constraints in the short term, increasing its constraint programming repertoire.

Another forthcoming and important application is running benchmarks to track the evolution of QUBO solvers and compare their performance against various problems.
This task will probably require specialized infrastructure to ensure consistency and reproducibility, for which the standards set by this work will be relevant.
In this sense, additional tools will be developed to archive and synthesize relevant QUBO instances for benchmarking.
Additionally, it is important to develop new packages that interact with solution-sampling devices to encompass an even broader range of devices and ensure diversity among solvers.
We envision that the ecosystem presented here will grow and allow an increasing number of operations research users to access quantum and quantum-inspired solution methods.

\section*{Acknowledgements}

The authors thank PSR for the great environment for research and open-source software development.
Special thanks to Mario Pereira and Sergio Granville for the motivation, support, and many discussions that led to this work.
We also thank Zachary Morrell for reading an early version of this manuscript and providing useful feedback.
DEBN acknowledges grant CCF 1918549 from the National Science Foundation (NSF) for funding.

\bibliographystyle{unsrtnat}
\bibliography{reference}

\appendix
% headings in the appendix read "Appendix A. Title", as in the T&F class
\titleformat{\section}[block]
  {\normalfont\large\bfseries}{Appendix~\thesection.}{0.5em}{}
% \section{Encoding methods and their coefficients}\label{sec:appendix:coefficients}

% Let \( f \in \mathscr{F} \) be of the form \eqref{eq:pbf}.
% The quantity \( \Delta(f) = \max_{\omega} \norm[1]{c_\omega} \), the absolute value of the coefficient with greater magnitude, is a rough yet interesting conditioning measure for the problem that $f$ represents.
% As a general recommendation, avoiding large values for $\Delta(f)$ leads to better reformulations.
% A large value in the largest coefficient affects the performance of the solvers as the precision required to represent the quadratic function accurately depends on the quotient of the smallest and largest coefficients in absolute value~(\cite{bian2014discrete}).

% This factor is one of the reasons why, in a situation where the unary encoding requires \( \bigo{n} \) variables to be implemented, it might still be preferred over the binary method, although it would require only \( \bigo{\log n} \) variables. 
{ 
\section{The Arithmetic Progression Encoding}\label{sec:appendix:ap}

Given an integer variable \( z \in [a, b] \), we can represent it as a linear combination of binary variables with an encoding based on the arithmetic progression.
Let \( n = b - a \) and \( N = \left\lceil{ \frac{1}{2} {\sqrt{1 + 8 n}} - \frac{1}{2} }\right\rceil\).
Now consider $N$ binary variables, \(\vec{x} \in \mathbb{B}^{N} \).
Finally, $z$ can be represented by the following linear combination of these binary variables:
\begin{equation}
    z =  a + \sum_{i = 1}^{N - 1} i\,x_{i} + \left({ n - \frac{N\,(N - 1)}{2}}\right)\,x_{N}
    \label{eq:ap}
\end{equation}
Note that the coefficients multiplying $x$ in the summation are in an arithmetic progression and thus bounded by \( N 
 = \bigo{\sqrt{n}} \) as \(n - \frac{1}{2} N (N - 1) \le  \frac{1}{2} \sqrt{1 + 8n} 
- \frac{5}{8} \le N \).
Naturally, it is possible to encode continuous variables by re-scaling \eqref{eq:ap}, $N$ fixed.

This strategy establishes a compromise solution between the resource-oriented binary encoding and the flat coefficients of the unary approach.
The bounded-coefficient technique also shares this kind of reasoning~(\cite{karimi2019practical}).

The Fixtars Amplify framework presents a similar proposal but with a slightly different formulation~(\cite{fixstars}).
To our knowledge, the two encoding schemes were developed independently, and no other presentation of an analogous method is found in the literature.

{

According to the documentation, for an interval of size $n$, their approach uses $N = 2\lfloor{\sqrt{n}}\rfloor + 1$ variables. We can prove that $N = \left\lceil{\frac{\sqrt{1 + 8n} - 1}{2}}\right\rceil$ is smaller.
}
}

\begin{figure}[h]
    \centering
    \includegraphics[width=0.8\linewidth]{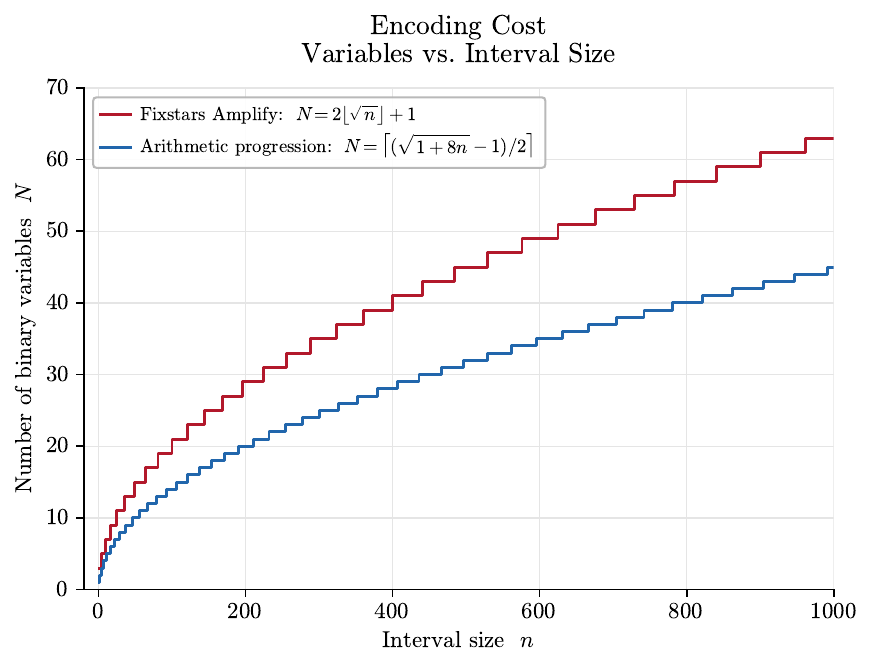}
    \caption{Encoding cost (number of binary variables $N$) versus interval size $n$
for the arithmetic-progression encoding of Eq.~\eqref{eq:ap} (blue) and the
Fixstars Amplify encoding (red). The former needs no more variables than the
latter for all $n$, with the gap growing as $n$ increases.}
    \label{fig:proof-sketch}
\end{figure}

\section{Benchmarking}\label{sec:appendix:benchmark}

The reformulation performance benchmark results presented in Figure \ref{fig:benchmark} were all obtained by single-threaded runs on the same computer.
% The details of its environment are summarized in Table \ref{tab:benchmark-specs}.
The computer was running on the \texttt{Linux Ubuntu 22.04} operating system, having \texttt{CPython 3.10.6} and \texttt{Julia 1.9.0} installed.
Its CPU is the model \texttt{Intel\textsuperscript{\textregistered} i7-7700K}, with base frequency \texttt{4.20 GHz}. 
Finally, the hardware has an installed memory capacity of \texttt{64 GB}.

Moreover, the set of software packages installed was \texttt{ToQUBO.jl}(v0.1.6), \texttt{PyQUBO}(v1.4.0), \texttt{OpenQAOA}(v0.1.3), \texttt{qubovert}(v1.2.5), \texttt{Qiskit}(v0.41.0), \texttt{amplify}(v0.11.1), \texttt{dimod}(v.0.12.14) and \texttt{docplex}(v.2.24.232).

All code used in the benchmarks, including visualization recipes, is available online in the \texttt{ToQUBO-benchmark} GitHub repository.\footnote{\url{https://github.com/JuliaQUBO/ToQUBO-benchmark}}

% \begin{center}
%     {\url{https://github.com/JuliaQUBO/ToQUBO-benchmark}}
% \end{center}

\end{document}